\magnification\magstep1
\newread\AUX\immediate\openin\AUX=\jobname.aux
\def\ref#1{\expandafter\edef\csname#1\endcsname}
\ifeof\AUX\immediate\write16{\jobname.aux gibt es nicht!}\else
\input \jobname.aux
\fi\immediate\closein\AUX
\def\today{\number\day.~\ifcase\month\or
  Januar\or Februar\or M\"arz\or April\or Mai\or Juni\or
  Juli\or August\or September\or Oktober\or November\or Dezember\fi
  \space\number\year}
\font\sevenex=cmex7
%\font\sevenex=cmex10 scaled 700
\scriptfont3=\sevenex
%\font\fiveex=cmex7 scaled 714
\font\fiveex=cmex10 scaled 500
\scriptscriptfont3=\fiveex
\def\epsilon{\varepsilon}
\def\theta{\vartheta}

\def\uauf{\lower1.7pt\hbox to 3pt{%
\vbox{\offinterlineskip
\hbox{\vbox to 8.5pt{\leaders\vrule width0.2pt\vfill}%
\kern-.3pt\hbox{\lams\char"76}\kern-0.3pt%
$\raise1pt\hbox{\lams\char"76}$}}\hfil}}
%%%%%%%%%%%%%%%%%%
% Makros f"ur Querverweise:
\def\cite#1{\expandafter\ifx\csname#1\endcsname\relax
{\bf?}\immediate\write16{#1 is not defined!}\else\csname#1\endcsname\fi}
\def\expandwrite#1#2{\edef\next{\write#1{#2}}\next}
\def\neverexpand{\noexpand\noexpand\noexpand}
\def\strip#1\ {}
\def\ncite#1{\expandafter\ifx\csname#1\endcsname\relax
{\bf?}\immediate\write16{#1 is not defined!}\else
\expandafter\expandafter\expandafter\strip\csname#1\endcsname\fi}
\newwrite\AUX
\immediate\openout\AUX=\jobname.aux
%%%%%%%%%%%%%%%%%%%%%
\newcount\Abschnitt\Abschnitt0
\def\beginsection#1. #2 \par{\advance\Abschnitt1%
\vskip0pt plus.10\vsize\penalty-250
\vskip0pt plus-.10\vsize\bigskip\vskip\parskip
\edef\TEST{\number\Abschnitt}
\expandafter\ifx\csname#1\endcsname\TEST\relax\else
\immediate\write16{#1 has already been defined!}\fi
\expandwrite\AUX{\neverexpand\ref{#1}{\TEST}}
\leftline{\bf\number\Abschnitt. \ignorespaces#2}%
\nobreak\smallskip\noindent\SATZ1}
%%%%%%%%%%%%%%%%%%
\def\Proof:{\par\noindent{\it Proof:}}
\def\Remark:{\ifdim\lastskip<\medskipamount\removelastskip\medskip\fi
\noindent{\bf Remark:}}
\def\Remarks:{\ifdim\lastskip<\medskipamount\removelastskip\medskip\fi
\noindent{\bf Remarks:}}
\def\Definition:{\ifdim\lastskip<\medskipamount\removelastskip\medskip\fi

\noindent{\bf Definition:}}
\def\Example:{\ifdim\lastskip<\medskipamount\removelastskip\medskip\fi
\noindent{\bf Example:}}
%%%%%%%%%%%%%%%%
\newcount\SATZ\SATZ1
\def\proclaim #1. #2\par{\ifdim\lastskip<\medskipamount\removelastskip
\medskip\fi
\noindent{\bf#1.\ }{\it#2}\Par
\ifdim\lastskip<\medskipamount\removelastskip\goodbreak\medskip\fi}
\def\Aussage#1{%
\expandafter\def\csname#1\endcsname##1.{\ifx?##1?\relax\else
\edef\TEST{#1\penalty10000\ \number\Abschnitt.\number\SATZ}
\expandafter\ifx\csname##1\endcsname\TEST\relax\else
\immediate\write16{##1 had already been defined!}\fi
\expandwrite\AUX{\neverexpand\ref{##1}{\TEST}}\fi
\proclaim {\number\Abschnitt.\number\SATZ. #1\global\advance\SATZ1}.}}
\Aussage{Theorem}
\Aussage{Proposition}
\Aussage{Corollary}
\Aussage{Conjecture}
\Aussage{Lemma}
%%%%%%%%%%%%%%%%
\font\la=lasy10
\def\strich{\hbox{$\vcenter{\hbox
to 1pt{\leaders\hrule height -0,2pt depth 0,6pt\hfil}}$}}
\def\dashedrightarrow{\hbox{%
\hbox to 0,5cm{\leaders\hbox to 2pt{\hfil\strich\hfil}\hfil}%
\kern-2pt\hbox{\la\char\string"29}}}

\def\Bindestrich{\penalty10000-\hskip0pt}
\let\_=\Bindestrich
\def\.{{\sfcode`.=1000.}}
%%%%%%%%%%%%%%%%%%%%%%%%%%%%%%%%%%%%

\def\Par{\par}
\def\:={\mathrel{\raise0,9pt\hbox{.}\kern-2,77779pt
\raise3pt\hbox{.}\kern-2,5pt=}}
\def\=:{\mathrel{=\kern-2,5pt\raise0,9pt\hbox{.}\kern-2,77779pt
\raise3pt\hbox{.}}}

\def\|#1|{\mathop{\rm#1}\nolimits}
\def\<{\langle}
\def\>{\rangle}
\let\Times=\times
\def\times{\mathop{\Times}}
\let\Otimes=\otimes
\def\otimes{\mathop{\Otimes}}
%%%%%%%%%%%%%%%%%%%%%%%%%%%%%%%%%
%Laden von Fonts:
\catcode`\@=11
\def\hex#1{\ifcase#1 0\or1\or2\or3\or4\or5\or6\or7\or8\or9\or A\or B\or
C\or D\or E\or F\else\message{Warnung: Setze hex#1=0}0\fi}
\def\fontdef#1:#2,#3,#4.{%
\alloc@8\fam\chardef\sixt@@n\FAM
\ifx!#2!\else\expandafter\font\csname text#1\endcsname=#2
\textfont\the\FAM=\csname text#1\endcsname\fi
\ifx!#3!\else\expandafter\font\csname script#1\endcsname=#3
\scriptfont\the\FAM=\csname script#1\endcsname\fi
\ifx!#4!\else\expandafter\font\csname scriptscript#1\endcsname=#4
\scriptscriptfont\the\FAM=\csname scriptscript#1\endcsname\fi
\expandafter\edef\csname #1\endcsname{\fam\the\FAM\csname text#1\endcsname}
\expandafter\edef\csname hex#1fam\endcsname{\hex\FAM}}
\catcode`\@=12 

%%%%%%%%%%%%%%%%%%%%%%%%%%%%%%%%%
\fontdef Ss:cmss10,,.
\fontdef Fr:eufm10,eufm7,eufm5.

			%Hier aufpassen!!!

%
\newread\AUXX
\immediate\openin\AUXX=msxym.tex
\ifeof\AUXX
\fontdef bbb:msbm10,msbm7,msbm5.
\fontdef mbf:cmmib10,cmmib7,.
\else
\fontdef bbb:msym10,msym7,msym5.
\fontdef mbf:cmmib10,cmmib10 scaled 700,.
\fi
\immediate\closein\AUXX

\def\FF{{\bbb F}}
\def\KK{{\bbb K}}
\def\QQ{{\bbb Q}}

\def\ZZ{{\bbb Z}}

\def\cH{{\cal H}}

\def\cP{{\cal P}}
\def\cS{{\cal S}}

\mathchardef\leer=\string"0\hexbbbfam3F
\mathchardef\subsetneq=\string"3\hexbbbfam24
\mathchardef\semidir=\string"2\hexbbbfam6E
\mathchardef\dirsemi=\string"2\hexbbbfam6F
\let\OL=\overline
\def\overline#1{{\hskip1pt\OL{\hskip-1pt#1\hskip-1pt}\hskip1pt}}

%<--                    Aufpassen  

%
%%%%%%%%%%%%
% Displayroutine
\abovedisplayskip 9.0pt plus 3.0pt minus 3.0pt
\belowdisplayskip 9.0pt plus 3.0pt minus 3.0pt
\newdimen\Grenze\Grenze2\parindent\advance\Grenze1em
\newdimen\Breite
\newbox\DpBox
\def\NewDisplay#1$${\Breite\hsize\advance\Breite-\hangindent
\setbox\DpBox=\hbox{\hskip2\parindent$\displaystyle{#1}$}%
\ifnum\predisplaysize<\Grenze\abovedisplayskip\abovedisplayshortskip
\belowdisplayskip\belowdisplayshortskip\fi
\global\futurelet\nexttok\WEITER}
\def\WEITER{\ifx\nexttok\qed\expandafter\leftQEDdisplay
\else\leftdisplay\fi}
\def\leftdisplay{\hskip-\hangindent\leftline{\box\DpBox}$$}
\def\leftQEDdisplay{\hskip-\hangindent
\line{\copy\DpBox\hfill\lower\dp\DpBox\copy\QEDbox}%
\belowdisplayskip0pt$$\bigskip\let\nexttok=}
\everydisplay{\NewDisplay}
%%%%%%%%%%%%
\newbox\QEDbox
\newbox\nichts\setbox\nichts=\vbox{}\wd\nichts=2mm\ht\nichts=2mm
\setbox\QEDbox=\hbox{\vrule\vbox{\hrule\copy\nichts\hrule}\vrule}
\def\qed{\leavevmode\unskip\hfil\null\nobreak\hfill\copy\QEDbox\medbreak}
%%%%%%%%%%%%%%
\newdimen\HIindent
\newbox\HIbox
\def\setHI#1{\setbox\HIbox=\hbox{#1}\HIindent=\wd\HIbox}
\def\HI#1{\par\hangindent\HIindent\hangafter=0\noindent\leavevmode
\llap{\hbox to\HIindent{#1\hfil}}\ignorespaces}
%%%%%%%%%%%%%%

\baselineskip12pt
\parskip2.5pt plus 1pt
\hyphenation{Hei-del-berg}
\def\L|Abk:#1|Sig:#2|Au:#3|Tit:#4|Zs:#5|Bd:#6|S:#7|J:#8||{%
\edef\TEST{[#2]}
\expandafter\ifx\csname#1\endcsname\TEST\relax\else
\immediate\write16{#1 hat sich geaendert!}\fi
\expandwrite\AUX{\neverexpand\ref{#1}{\TEST}}
\HI{[#2]}
\ifx-#3\relax\else{#3}: \fi
\ifx-#4\relax\else{#4}{\sfcode`.=3000.} \fi
\ifx-#5\relax\else{\it #5\/} \fi
\ifx-#6\relax\else{\bf #6} \fi
\ifx-#8\relax\else({#8})\fi
\ifx-#7\relax\else, {#7}\fi\Par}

\def\B|Abk:#1|Sig:#2|Au:#3|Tit:#4|Reihe:#5|Verlag:#6|Ort:#7|J:#8||{%
\edef\TEST{[#2]}
\expandafter\ifx\csname#1\endcsname\TEST\relax\else
\immediate\write16{#1 hat sich geaendert!}\fi
\expandwrite\AUX{\neverexpand\ref{#1}{\TEST}}
\HI{[#2]}
\ifx-#3\relax\else{#3}: \fi
\ifx-#4\relax\else{#4}{\sfcode`.=3000.} \fi
\ifx-#5\relax\else{(#5)} \fi
\ifx-#7\relax\else{#7:} \fi
\ifx-#6\relax\else{#6}\fi
\ifx-#8\relax\else{ #8}\fi\Par}

\def\Pr|Abk:#1|Sig:#2|Au:#3|Artikel:#4|Titel:#5|Hgr:#6|Reihe:{%
\edef\TEST{[#2]}
\expandafter\ifx\csname#1\endcsname\TEST\relax\else
\immediate\write16{#1 hat sich geaendert!}\fi
\expandwrite\AUX{\neverexpand\ref{#1}{\TEST}}
\HI{[#2]}
\ifx-#3\relax\else{#3}: \fi
\ifx-#4\relax\else{#4}{\sfcode`.=3000.} \fi
\ifx-#5\relax\else{In: \it #5}. \fi
\ifx-#6\relax\else{(#6)} \fi\PrII}
\def\PrII#1|Bd:#2|Verlag:#3|Ort:#4|S:#5|J:#6||{%
\ifx-#1\relax\else{#1} \fi
\ifx-#2\relax\else{\bf #2}, \fi
\ifx-#4\relax\else{#4:} \fi
\ifx-#3\relax\else{#3} \fi
\ifx-#6\relax\else{#6}\fi
\ifx-#5\relax\else{, #5}\fi\Par}
\setHI{[KKLV]\ }
\sfcode`.=1000

\def\alp{\alpha} 
\def\bet{\beta} 
\def\lam{\lambda}
\def\sig{\sigma}
\def\gam{\gamma}

\def\abar{{\overline\alp}}
\def\bbar{{\overline\bet}}

\def\lbar{{\overline\lam}}
\def\mbar{{\overline\mu}}
\def\vbar{\overline v}

\def\atil{{\widetilde\alp}}
\def\btil{{\widetilde\bet}}
\def\gtil{{\widetilde\gam}}
\def\mtil{{\widetilde\mu}}

\def\qt{\left[\alp\atop\bet\right]_{q,t}}
\def\abr{\left[\alp\atop\bet\right]_r}
\def\qtinv{\left[\alp\atop\bet\right]_{1/q,1/t}}
\def\lmr{\left(\lam\atop\mu\right)_r}

\fontdef Ss:cmss10,,.
\font\BF=cmbx10 scaled \magstep1
\font\CSC=cmcsc10 %scaled \magstephalf
%\ignore

\baselineskip15pt
%\ignore

{\baselineskip1.5\baselineskip\rightskip0pt plus 5truecm
\leavevmode\vskip0truecm\noindent
\BF The binomial formula for nonsymmetric Macdonald polynomials
}
\vskip1truecm
\leftline{{\CSC Siddhartha Sahi}%
\footnote*{\rm This work was supported by an NSF grant}}
\leftline{Department of Mathematics, Rutgers University, 
New Brunswick NJ 08903, USA
}
\bigskip
\beginsection intro. Introduction

The $q$-binomial theorem is essentially the expansion of 
$(x-1)(x-q) \cdots(x-q^{k-1})$ in terms of the monomials $x^d$. 
In a recent paper \cite{O}, A.~Okounkov has proved a beautiful 
multivariate generalization of this in the context of symmetric 
Macdonald polynomials \cite{M1}. These polynomials have 
nonsymmetric counterparts \cite{M2} which are of substantial 
interest, and in this paper we establish nonsymmetric analogues 
of Okounkov's results. 

An integral vector $v\in\ZZ^n$ is called ``dominant'' if 
$v_1\ge\cdots\ge v_n$; and it is called a ``composition'' 
if $v_i\ge0$, for all $i$. To avoid ambiguity we reserve the 
letters $u,v$ for integral vectors, $\alp,\bet,\gam$ for 
compositions, and $\lam,\mu$ for ``partitions'' (dominant 
compositions).

We write $|v|$ for $v_1+\cdots+v_n$, and denote by $w_v$ the 
(unique) shortest permutation in the symmetric group $S_n$ 
such that $v^+= w_v^{-1}(v)$ is dominant. Let $\FF$ be field 
$\QQ(q,t)$ where $q,t$ are indeterminates. 
We write $\tau=(1,t^{-1},\cdots,t^{-n+1})$, and define $\vbar=
\vbar(q,t)$ in $\FF^n$ by 
$$
\vbar_i= q^{v_i}(w_v\tau)_i.
$$

Inhomogeneous analogues of nonsymmetric 
Macdonald polynomials were introduced in \cite{K} and \cite{S2}.
They form an $\FF$-basis for $\FF[x]=\FF[x_1,\cdots,x_n]$ and
% for the subspace of symmetric polynomials, respectively, 
 are defined as follows:

\Definition: $G_\alp\equiv G_\alp(x;q,t)$ is the unique 
polynomial of degree $\le |\alp|$ in $\FF[x]$ such that
\item{1)} the coefficient of $x^\alp\equiv x_1^{\alp_1}\cdots 
x_n^{\alp_n}$ in $G_\alp$ is $1$,
\item{2)} $G_\alp$ vanishes at $x=\bbar$, for all 
compositions $\bet\ne\alp$ such that $|\bet|\le |\alp|$.

As shown in Th.~3.9 of \cite{K}, the top homogeneous part of $G_\alp$
is the nonsymmetric Macdonald 
polynomial $E_\alp$ for the root system $A_{n-1}$ (\cite{M2} and \cite{C}). 
Moreover by Th.~4.5 of \cite{K} we have $G_\alp(\bbar)=0$ unless 
``$\alp\subseteq\bet$''. Here $\alp\subseteq\bet$ means that if we 
write $w=w_\bet w_\alp^{-1}$ then $\alp_i<\bet_{w(i)}$ if $i<w(i)$ 
and $\alp_i\le\bet_{w(i)}$ if $i\ge w(i)$.

In this paper we obtain several new results about
the polynomials $G_\alp$. 

Our first result is a formula for the special value $G_\alp(a\overline0)=
G_\alp(a\tau)\in\FF[a]$ where $a$ is an indeterminate. This can be 
described in the following manner: 

We identify $\alp$ with the ``diagram'' consisting of 
points $(i,j)\in\ZZ^2$ with $1\le i\le n$ 
and $1\le j\le\alp_i$. For $s=(i,j)\in\alp$ we define the 
{\it arm\/},  {\it leg\/}, {\it coarm\/}, and {\it coleg\/} 
of $s$ by

$a(s)=\alp_i-j,\;
l(s)= \#\{k>i\mid j\le\alp_k\le\alp_i\}+
\#\{k<i\mid j\le\alp_k+1\le\alp_i\},$

$a'(s)=j-1,\;l'(s)= \#\{k>i\mid \alp_k>\alp_i\}+
\#\{k<i\mid \alp_k\ge\alp_i\}.$

\Theorem eval.  $\displaystyle G_\alp(a\tau)= 
\prod_{s\in\alp} 
\left({t^{1-n}-q^{a'(s)+1}t^{1-l'(s)} 
\over
 1-q^{a(s)+1}t^{l(s)+1}}\right)
\prod_{s\in\alp} \left(at^{l'(s)}-q^{a'(s)}\right).$

Let $w_o$ be the longest element of $S_n$ (which 
interchanges each $i$ with $n-i+1$), and put $\btil= 
\overline{-w_o\bet}$ and $\bbar^{-1}=\bbar(q^{-1},t^{-1})=
(\bbar_1^{-1},\cdots,\bbar_n^{-1})$. Then we have the following 
crucial ``reciprocity'' result:

\Theorem Oko. There is a (unique) polynomial $O_\alp$ of degree 
$\le|\alp|$ in $\QQ(q,t,a)[x]$ such that 
$O_\alp(\bbar^{-1})=G_\bet(a\atil)/G_\bet(a\tau)$ 
for all $\bet$. 

We now introduce the following variants of $G_\alp$, which also form 
a basis for $\FF[x]$:

\Definition: $G'_\alp= G'_\alp(x;q,t)$ is the 
unique polynomial in $\FF[x]$ such that
\item{1)} $G'_\alp$ and $G_\alp$ have the same top degree terms, i.e. $E_\alp$, 
\item{2)} $G'_\alp$ vanishes at $x=\btil$ for all $\bet$ with $|\bet|<|\alp|$.

The existence of $G'_\alp$ can be proved along the same lines
as that of $G_\alp$ (\cite{K} Th. 2.3, \cite{S2} Th. 4.3). 
One verifies that polynomials of degree $\le d$ are 
uniquely determined by their values at $x=\btil$
for $|\bet|\le d$. Hence the lower degree terms
of $G'_\alp$ are determined by 2).

\Definition: The ``nonsymmetric
$(q,t)$-binomial coefficients'' are defined by $$\qt
= {G_\bet(\abar)\over G_\bet(\bbar)}\equiv {G_\bet(\abar(q,t);q,t)
\over G_\bet(\bbar(q,t);q,t)}.$$

Our main result is the following relationship between $G_\alp$ 
and $G'_\bet$:

\Theorem binom. 
$\displaystyle {G_\alp(ax)\over G_\alp(a\tau)} = 
\sum_{\bet\subseteq\alp} 
a^{|\bet|}\qtinv
{G'_\bet(x)\over G_\bet(a\tau)}.$ 

\Corollary first. $\displaystyle{ G_\alp(x)\over G_\alp(0)}
= \sum_{\bet\subseteq\alp} \qtinv
{E_\bet(x)\over G_\bet(0)}.$ \qed

\Corollary Gprime.
$\displaystyle {E_\alp(x)\over E_\alp(\tau)}
= \sum_{\bet\subseteq\alp} \qtinv
{G'_\bet(x)\over E_\bet(\tau)}.$ \qed

The corollaries follow from \cite{binom} by (1) replacing $x$ by 
$a^{-1}x$ and letting $a \rightarrow0$, and (2) by letting $a\rightarrow\infty$.
For $n=1$, \cite{first} is essentially the $q$-binomial theorem.

%which is independent of $t$. Thus $\displaystyle 
%\left[d\atop k \right]_{q,t}= 
%\prod_{i=0}^{k-1} {q^d-q^i\over q^k-q^i} =\prod_{i=0}^{k-1} 
%{1-q^{d-i}\over 1-q^{k-i}}$ is the usual $q$-binomial coefficient,

If we put $t=q^r$ and let $q\rightarrow1$ then
$E_\alp(x;r)\equiv\lim_{q\rightarrow1} E_\alp(x;q,q^r)$ 
is the nonsymmetric Jack polynomial \cite{Op}.
To discuss this limiting case, we define 
$\delta\equiv(0,-1,\cdots,-n+1)$, $\rho=r\delta$, and
$\abar(r)=\alp+w_\alp\rho$.

\Definition: $G_\alp(x;r)$ is the unique 
polynomial of degree $\le |\alp|$ in $\QQ(r)[x]$ such that
\item{1)} the coefficient of $x^\alp\equiv x_1^{\alp_1}\cdots 
x_n^{\alp_n}$ in $G_\alp(x;r)$ is $1$,
\item{2)} $G_\alp(x;r)$ vanishes at $x=\bbar(r)$, for all 
compositions $\bet\ne\alp$ such that $|\bet|\le |\alp|$.

\Definition: The ``nonsymmetric $r$-binomial coefficients'' are 
$\displaystyle\abr = {G_\bet(\abar(r);r)\over G_\bet(\bbar(r);r)}.$

\Definition: $G'_\alp(x;r)$ is the 
unique polynomial in $\QQ(r)[x]$ such that
\item{1)} $G'_\alp(x;r)$ and $G_\alp(x;r)$ have the same top degree terms,
\item{2)} $G'_\alp(x;r)$ vanishes at $x=\btil(r)\equiv \overline{-w_o\bet}(r)$ 
for all $\bet$ with $|\bet|<|\alp|$.

Theorems \ncite{eval} --- \ncite{binom} have analogues in this setting.

If $a$ is a scalar and $x$ is a vector, write $a+x$
for $(a+x_1,\cdots,a+x_n)$.  

\Theorem eval2. $\displaystyle G_\alp(a+\rho;r)= 
\prod_{s\in\alp}\left({a'(s)+1-rl'(s)+rn \over
a(s)+1+rl(s)+r}\right)
\prod_{s\in\alp}\left(a-a'(s)+rl'(s)\right)$.

\Theorem Oko2. There is a (unique) polynomial $O_\alp(x;r)$ of 
degree $\le |\alp|$ in $\QQ(a,r)[x]$ such that we have 
$\displaystyle O_\alp(\bbar(r);r)= G_\bet(a+\atil(r);r)/
G_\bet(a+\rho;r)$ for all $\bet$. 

\Theorem binom2.
$\displaystyle {G_\alp(a+x;r)\over G_\alp(a+\rho;r)}
= \sum_{\bet\subseteq\alp} \abr 
{G'_\bet(x;r)\over G_\bet(a+\rho;r)}.$ 

Since $\abar_i(r)= \lim_{q\rightarrow1}(\abar_i(q,q^r)-1)/(q-1)$, 
as in \cite{K} Th. 6.2. we get $G_\alp(x;r)\equiv\lim_{q\rightarrow1} 
G_\alp(1+(q-1)x;q,q^r)/(q-1)^{|\alp|}$.
It follows that the top terms of $G_\alp(x;r)$ and $G'_\alp(x;r)$ are  
$E_\alp(x;r)$, and that
$\abr= \lim_{q\rightarrow1} \left[\alp \atop\bet\right]_{q,q^r}
= \lim_{q\rightarrow1} \left[\alp \atop\bet\right]_{1/q,1/q^r}$

So setting $x=ax$ and letting $a\rightarrow\infty$ in \cite{binom2} we get

\Corollary Las. 
$\displaystyle {E_\alp(1+x;r)\over E_\alp(1;r)}
= \sum_{\bet\subseteq\alp} \abr 
{E_\bet(x;r)\over E_\bet(1;r)}.$ \qed

It seems to be difficult to deduce Theorems \ncite{eval2} --- 
\ncite{binom2} directly from Theorems \ncite{eval} --- \ncite{binom} 
by a limiting procedure. However the {\it proofs\/} of in the 
$(q,t)$-case {\it can\/} be modified to make them work in this 
setting.

We now describe some {\it new\/} phenomena in the limiting case.
Write $s_i$ for the transposition $(i\;i+1)$ which acts
on $\QQ(a)[x]$ by permuting $x_i$ and $x_{i+1}$, and let  
$$\sig_i=s_i+{r\over x_i-x_{i+1}}(1-s_i).$$
Then, as observed in \cite{K} Cor. 6.5, the map $\sig:s_i\mapsto 
\sig_i$ extends to a representation of $S_n$. 

\Theorem relate. 
$\displaystyle G'_\alp(x;r)=(-1)^{|\alp|}
\sig(w_o)w_oG_\alp(-x-(n-1)r;r).$

Using this, and writing $G^+_\alp(x;r):= (-1)^{|\alp|}
G_\alp(-x-(n-1)r;r)$, we get

\Corollary plus.
$\displaystyle {\sig(w_o)G_\alp(a+x;r)\over G_\alp(a+\rho;r)}
= \sum_{\bet\subseteq\alp} \abr
{w_oG^+_\bet(x;r)\over G_\bet(a+\rho;r)}.$ \qed

As mentioned earlier, the symmetric analogues of Theorems
\ncite{eval} --- \ncite{binom} have been established in \cite{O}.
In the case of symmetric Jack polynomials, expansions in the form 
of \cite{Las} were first considered by \cite{Bi} ($r=1/2$), 
and \cite{Lc} ($r=1$), and in general by \cite{Ls}. The analogues 
of Theorems \ncite{eval2}, \ncite{Oko2}, and \cite{Las} have 
been obtained by \cite{OO}, but the analogues of \cite{binom2}, 
and \ncite{relate} seem not to have been considered by them. Since 
these follow easily by our techniques, we shall formulate and prove 
them in \cite{binom3} and \cite{relate2} below.

While our proof follows the same general outline as Okounkov's 
argument, there are several differences. First, a decisive role 
is played by the affine Hecke algebra and Cherednik operators, 
and the Hecke recursions satisfied by the $G_\alp$ actually yield 
a simplification of part of the argument. On the other hand, there 
are some subtleties in the nonsymmetric case, as exhibited by 
the definition of $G_\alp'$.

\beginsection prelim. Preliminaries

We start by recalling certain basic properties of the 
$G_\alp(x;q,t)$ (see \cite{K} and \cite{S2}).

The main result of \cite{K} (Th.~3.6) is that the $G_\alp$ 
satisfy the eigen-equations  
$$\Xi_i G_\alp= \abar_i^{-1}G_\alp$$ 
for the ``inhomogeneous Cherednik operators'' defined
by $$\Xi_i= x_i^{-1}+x_i^{-1}H_i\ldots H_{n-1}\Phi H_1\ldots 
H_{i-1}.$$
In turn, the operators $\Phi$ and $H_i$ are defined by
$$\eqalign{\Phi f(x_1,\cdots,x_n)&=(x_n-t^{-n+1})
                 f(x_n/q,x_1, \cdots,x_{n-1})\cr
H_i&=ts_i-(1-t){x_i\over x_i-x_{i+1}}(1-s_i).}$$ 

The $H_i$'s satisfy the braid relations and the 
identity $(H_i-t)(H_i+1)=0$, and generate a representation 
of the Iwahori-Hecke algebra $\cH$ of $S_n$ on $\FF[x]$.

Next, write $v^\#=(v_n-1,v_1,\cdots,v_{n-1})$; and let $a$ 
be an indeterminate. 

\Lemma discr.
\item{$1)$} $\Phi f(a\vbar)= (a\vbar_n-t^{-n+1})f(a\overline{v^\#})$
\item{$2)$} $\displaystyle H_if(a\vbar)= {(t-1)\vbar_i \over \vbar_i- \vbar_{i+1}}f(a\vbar)
+{\vbar_i -t\vbar_{i+1} \over \vbar_i- \vbar_{i+1}}
f(a\overline{s_iv}).$ 

This is proved just as in \cite{K} Lemmas 2.1, 3.1. The main 
point in 2) is that for $v\in\ZZ^n$, $s_iv=v \Rightarrow 
\vbar_i-t\vbar_{i+1}=0$, and $s_iv\ne v \Rightarrow s_i\vbar=
\overline{s_i v}$. 

\Lemma recur. 
\item {$1)$} If $\alp_n>0$ then $G_\alp=q^{\alp_n-1} \Phi G_{\alp^\#}$.
\item {$2)$} If $\alp_i> \alp_{i+1}$ then $G_\alp= 
(H_i+ (1-t)d^{-1})G_{s_i\alp}$ where $d= (1-\abar_i/\abar_{i+1})$. 

This is essentially in \cite{K} and \cite{S2} --- here is a sketch of the 
argument: Evidently the right sides of 1) and 2) have degree $\le|\alp|$,
and by using \cite{discr} one verifies the vanishing conditions. It remains 
only to check that the coefficient of $x^\alp$ is $1$. This obvious for 1), 
while for 2) one has to use the triangularity of $\Xi_i$ (Lemma 3.10 of 
\cite{K}). 

In connection with \cite{eval}, we define scalars $d_\alp(q,t)=
\prod_{s\in\alp} (1-q^{a(s)+1}t^{l(s)+1})$, $e_\alp(q,t)=\prod_{s\in\alp} 
 (t^{1-n}-q^{a'(s)+1}t^{1-l'(s)})$, and 
$\phi_\alp(a;q,t)= \prod_{s\in\alp} \left(at^{l'(s)}-q^{a'(s)}\right).$

\Lemma derecur.
\item {$1)$} If $\alp_n>0$ then $d_\alp/d_{\alp^\#}= 1-t^n\abar_n$, 
$e_\alp/e_{\alp^\#}=t^{1-n}-t\abar_n$, $\phi_\alp(0)=
-q^{\alp_n-1}\phi_{\alp^\#}(0)$.
\item {$2)$} If $\alp_i> \alp_{i+1}$ then $d_\alp= 
{1- \abar_i/\abar_{i+1}\over1-t\abar_i/\abar_{i+1}}d_{s_i\alp}$.
\item {$3)$}  $e_{w\alp}=e_\alp$ and 
$\phi_{w\alp}=\phi_\alp$ for all $w$ in $S_n$.

The lemma can be proved in a manner very similar to
Lemmas 4.1 and 4.2 in \cite{S3}. To illustrate the argument,
we sketch the proof of $e_\alp/e_{\alp^\#}=t^{1-n}-t\abar_n$,
other proofs are similar: 
It follows from the definition of $\abar$ that $\abar_i=
q^{\alp_i}t^{-k_i}$ where $k_i= \#\{k<i\mid\alp_k\ge\alp_i\}+
\#\{k>i\mid\alp_k>\alp_i\}.$

The diagram of $\alp$ is obtained from $\alp^\#$
by adding a point to the end of the first row, and moving
this row to the last place. The new point $s=(n,\alp_n)\in\alp$ 
has $a'(s)=\alp_n-1$ and $l'(s)=\#\{k<n\mid\alp_k>\alp_n\}=k_n$,
while coarms and colegs of other points are unchanged.
Thus \quad $e_\alp/e_{\alp^\#}= 
t^{1-n}-q^{a'(s)+1}t^{1-l'(s)} = t^{1-n}-q^{\alp_n}t^{1-k_n}=
t^{1-n}-t\abar_n.$

We shall also need limit versions of these results 
which are proved similarly:

First, by \cite{K} Th. 6.6, we know that the $G_\alp(x;r)$ satisfy 
the eigen-equations $$\widetilde \Xi_iG_\alp(x;r)= \abar(r)G_\alp,$$ 
where the ``limit'' Cherednik operators are defined by 
$$\widetilde 
\Xi_i:= x_i-\sig_i\ldots \sig_{n-1}\widetilde\Phi \sig_1\ldots
\sig_{i-1}.$$
where $\sig_i=s_i+r(x_i-x_{i+1})^{-1}(1-s_i)$ is as in the previous 
section, and 
$$ \widetilde\Phi f(x)=(x_n+(n-1)r)f(x_n-1,x_1\ldots,x_{n-1}).
$$

\Lemma discr2.
\item{$1)$} $\widetilde\Phi f(a+\vbar)= 
(a+\vbar_n+nr-r)f(a+\overline{v^\#})$
\item{$2)$} $\displaystyle \sig_if(a+\vbar)= 
{r \over \vbar_i- \vbar_{i+1}}f(a+\vbar)
+{\vbar_i -\vbar_{i+1}-r \over \vbar_i- \vbar_{i+1}}
f(a+\overline{s_iv}).$ \qed

\Lemma recur2.
\item {$1)$} If $\alp_n>0$ then $G_\alp= \widetilde\Phi G_{\alp^\#}$.
\item {$2)$} If $\alp_i> \alp_{i+1}$ then $G_\alp= 
(\sig_i+rd^{-1})G_{s_i\alp}$ where $d=\abar_i-\abar_{i+1}$. \qed

In connection with \cite{eval2} we define scalars $d_\alp(r)=
\prod_{s\in\alp}\left( a(s)+1+rl(s)+r\right)$,
$e_\alp(r)=\prod_{s\in\alp}\left(a'(s)+1-rl'(s)+rn\right)$,
and $\phi_\alp(a;r)=\prod_{s\in\alp}\left(a-a'(s)+rl'(s)\right)$.

\Lemma derecur2.
\item {$1)$} If $\alp_n>0$ then $d_\alp(r)/d_{\alp^\#}(r)= 
rn+\abar_n(r) =e_ \alp(r)/e_{\alp^\#}(r)$.
\item {$2)$} If $\alp_i> \alp_{i+1}$ then $d_\alp(r)= 
{ d\over d+r}d_{s_i\alp}(r)$, where $d=\abar_i-\abar_{i+1}$. 
\item {$3)$}  $e_{w\alp}(r)=e_\alp(r)$ and 
$\phi_{w\alp}=\phi_\alp$ for all $w$ in $S_n$.  \qed

We now briefly discuss the symmetric case. 

\Definition: $R_\lam(x;q,t)$ is the unique symmetric polynomial 
of degree $\le |\lam|$ which vanishes at $x=\mbar$ for partitions 
$\mu\ne\lam, |\mu|\le|\lam|$ and is normalized so that the 
coefficient of $x^\lam$ is $1$. 

\Definition: $R_\lam(x;r)$ is the unique symmetric polynomial 
of degree $\le |\lam|$ which vanishes at $x=\mbar(r)$ for partitions 
$\mu\ne\lam, |\mu|\le|\lam|$ and is normalized so that the 
coefficient of $x^\lam$ is $1$. 

The existence and uniqueness of $R_\lam(x;q,t)$ was proved 
in \cite{K} and \cite{S2}, as was the fact that its top term
is the Macdonald polynomial $P_\lam(q,t)$. In the case of
$R_\lam(x;r)$ these results were established in \cite{S1}
and \cite{KS}.

As in \cite{S2} Th. 4.6 and \cite{K} Cor. 2.6, we have:
\Lemma RG. Let $V_\lam$ be the $\FF$-span of
$\{E_\alp(x;q,t)\mid\alp^+=\lam\}$. Then $V_\lam$
is a module for the Hecke algebra $\cH$, and
$V_\lam^\cH=\FF R_\lam(x;q,t)$. 

\Lemma RG2. Let $V_\lam(r)$ be the $\QQ(r)$-span of
$\{E_\alp(x;r)\mid\alp^+=\lam\}$. Then $V_\lam(r)$
is a module for $\sig(S_n)$, and
$V_\lam(r)^{\sig(S_n)}=\QQ(r) R_\lam(x;r)$. 

Finally, For compatibility of notation between \cite{K}, 
\cite{O} and \cite{S2} we point out that 

\item{1)} \cite{K} uses $P_\lam$ for $R_\lam$, $\overline 
P_\lam$ for $P_\lam$, $\overline{E}_\alp$ for $E_\alp$, and 
$E_\alp$ for $G_\alp$.

\item{2)} \cite{O} uses $P_\lam^*(x)$ for the ``(shifted)'' 
polynomial $R_\lam(x\tau)\equiv R_\lam(x_1, x_2t^{-1}, \cdots, 
x_nt^{1-n})$ which vanishes at $(q^\mu_1,\cdots,q^\mu_n)$ 
and is symmetric in the variables $x_it^{-i}$.

\item{3)} \cite{S2} uses $R_\lam(x;q,t)$ to denote the polynomial
 $t^{-(n-1)|\lam|}R_\lam(xt^{n-1};q^{-1},t^{-1})$, which is 
symmetric and vanishes at the points $x=(q^{-\mu_1}t^{-n+1},
\cdots,q^{-\mu_{n-1}}t^{-1}, q^{-\mu_n})$. Its top term is 
$P_\lam(x;q^{-1},t^{-1})$ which equals $P_\lam(x;q,t)$ by \cite{M1}.

\beginsection Eval. Evaluation

In this section we prove the evaluation formulas
\cite{eval} and \cite{eval2}.

\Lemma inva. For all $w\in S_n$, we have $d_{w\alp}(q,t)G_{w\alp}(a\tau;q,t) 
=d_\alp(q,t)G_\alp(a\tau;q,t).$ 

\Proof: It suffices to verify this for $w=s_i$, and we may also assume 
that $\alp_i>\alp_{i+1}$.

 Since $\tau=\overline 0$, substituting $v=0$ in 2) 
of \cite{discr} we get $(H_if)(a\tau)=tf(a\tau)$ for all functions $f$.
Combining this with 2) of \cite{recur} we get
$$ G_\alp(a\tau)= (t+ (1-t)d^{-1})G_{s_i\alp}(a\tau)
={1-t \abar_i/\abar_{i+1}\over1-\abar_i/\abar_{i+1}}
G_{s_i\alp}(a\tau).$$ 
The result now follows from part 2) of \cite{derecur}. \qed

\cite{eval} states $\displaystyle 
d_\alp G_\alp(a\tau)= e_\alp\phi_\alp(a\tau)$, 
and we first establish this for $a=0$.

\Lemma zerosp. $\displaystyle d_\alp G_\alp(0)= e_\alp\phi_\alp(0)=
e_\alp \prod_{s\in\alp}\bigl(-q^{a'(s)}\bigr). $

\Proof: The case $\alp=0$ is trivial, and we proceed
by induction on $|\alp|$ assuming $\alp\ne0$. 
By \cite{inva} and part 3) of \cite{derecur} both sides are 
unchanged if permute $\alp$, so we may assume that 
$\alp_n>0$ and that $d_{\alp^\#} G_{\alp^\#}(0)=e_{\alp^\#}
\phi_{\alp^\#}(0)$. Thus it suffices to prove 
$$
{G_\alp(0)\over G_{\alp^\#}(0)}= 
\left({e_\alp \over e_{\alp^\#}}\right)
\left({d_{\alp^\#} \over d_\alp }\right)
\left({\phi_\alp(0) \over \phi_{\alp^\#}(0)}\right).$$

The left side can be computed by combining 1) of 
Lemmas \ncite{discr} and \ncite{recur}, and
the right side can be computed by 1) of \cite{derecur}.
In each case we get $-q^{\alp_n-1} t^{1-n}$.  \qed

We now deduce \cite{eval} from the symmetric case \cite{O}.

\Proof: (of \cite{eval}) If $\lam$ is a partition 
then, by \cite{O} (formula (1.9)), $R_\lam(a\tau)$ is an 
$\FF$-multiple of $\phi_\lam(a)$. 

Next, if $\alp$ be a composition such that $\alp^+=\lam$, then
by \cite{RG} there are some coefficients $c_w\in\FF$ 
such that $R_\lam(x)=\sum_{w\in S_n}c_w d_{w\alp} G_{w\alp}(x)$. 
Evaluating at $x=a\tau$ and using \cite{inva} we get
$R_\lam(a\tau)=(\sum c_w)d_\alp G_\alp(a\tau).$
It follows that $d_\alp(q,t)G_\alp(a\tau)$ is an $\FF$-multiple of 
$\phi_\lam(a)=\phi_\alp(a)$. 

Setting $a=0$ and using \cite{zerosp} we see that this multiple
is $e_\alp(q,t)$ and \cite{eval} follows. \qed

\Proof: (of \cite{eval2}) Arguing as in \cite{inva} we deduce
that $d_{w\alp}(r)G_{w\alp}(a+\rho;r) =d_\alp(r)G_\alp(a+\rho;r)$.
Next, by formula 2.3 of \cite{OO}, $R_\lam(a+\rho;r)$ is a $\QQ(r)$-
multiple of $\phi_\lam(a;r)$. Arguing as before, we conclude that
$d_\alp(r) G_\alp(a+\rho;r)$ is a $\QQ(r)$-multiple of $\phi_\alp(a;r)$.
 
Letting $a\rightarrow\infty$ we see that the multiple
is $d_\alp(r) E_\alp(1;r)$  which equals $e_\alp(r)$ 
by Th. 1.3 of \cite{S3}. The result follows.  \qed

\beginsection recip. Reciprocity 

In this section we prove \cite{Oko} and \cite{Oko2}. 

\Proof: (of \cite{Oko}) Write $\KK=\FF(a)=\QQ(q,t,a)$.
For $f$ in $\KK[x]$, we have
 $$\Xi_if(a\vbar)=(a\vbar_i)^{-1}f(a\vbar)+
(a\vbar_i)^{-1} H_i\ldots H_{n-1}\Phi H_1\ldots H_{i-1}f(a\vbar).$$
Since $|v^\#|=|v|-1$ and $|s_iv|=|v|$, it follows from \cite{discr}
that the second term on the right is a combination of 
$f(a\overline u)$ with $|u|=|v|-1$, where the 
coefficients do not depend on $f$. 
Thus if $p$ is a polynomial of degree $d$, and we write
$p(\Xi)\equiv p(\Xi_1,\ldots,\Xi_n)$, then
$p(\Xi)f(a\tau)\equiv p(\Xi)f(a\overline 0)=
\sum_{|\bet|\le d}c_p(\bet) f(a\btil)$,
with coefficients $c_p(\bet)$ independent of $f$.

Let $\cP$ be the space of polynomials in $\KK[x]$ of degree
$\le d$ and let $\cS$ be the set of compositions $\bet$
in $\ZZ_+^n$ with $|\bet|\le d$. Then $p\mapsto c_p$ is
a $\KK$-linear map from $\cP$ to $\KK^{\cS}$, and we claim that 
this map is bijective. 

Since the spaces have the same dimension, it suffices to 
check injectivity. If $c_p=0$ then 
$p(\Xi)f(a\tau)=0$ for all $f$. In particular, 
setting $f=G_\bet$ we obtain $p(\bbar^{-1})G_\bet(a\tau)=0$. 
By \cite{eval} $G_\bet(a\tau)\ne0$, and  
it follows that $p$ vanishes at the points $\bbar^{-1}= 
\bbar(q^{-1},t^{-1})$ for
all $\bet$, and hence $p=0$, proving injectivity. 

Now fix $\alp$ with $|\alp|=d$, and  let $O_\alp$ be the polynomial 
in $\cP$ whose image under $p\mapsto c_p$ is the delta function at 
$\alp$ in $\KK^{\cS}$.  Then $O_\alp$ has degree $\le |\alp|$, and
satisfies $O_\alp(\Xi)f(a\tau)=f(a\atil)$ for all $f$. Setting 
$f=G_\bet$ we get $O_\alp( \bbar^{-1})G_\bet(a\tau)
=G_\bet(a\atil)$. \qed

\Proof: (Of \cite{Oko2}) This is proved similarly by using
the limit Cherednik operators $\widetilde \Xi_i$ and \cite{discr2}. 
\qed

\beginsection Main. The binomial formula

We now prove \cite{binom} and \cite{binom2}.

\Proof: (of \cite{binom}) Since the $G'_\bet$ form a basis for 
$\KK[x]$, there exist $b_{\bet\alp}\in\KK$ such that
$$ (*) \qquad {G_\alp(ax)\over G_\alp(a\tau)}
=\sum_{\bet:|\bet|\le|\alp|} b_{\bet\alp} G'_\bet(x).$$ 
Substituting $x=\gtil$ and using \cite{Oko} we get 
$O_\gam(\abar^{-1})=\sum_\bet b_{\bet\alp}G'_\bet(\gtil)$.

Let $G$ be the (infinite) matrix whose entries are
$g_{\gam\bet}=G'_\bet(\gtil)$. By Th.~4.2
in \cite{S2} polynomials of degree $\le d$ 
are determined by their values at the points 
$\{\gtil:|\gam|\le d\}$, and it follows that 
$G$ has an inverse $H$, and we get $b_{\bet\alp}=
\sum_\gam h_{\bet\gam} O_\gam(\abar^{-1}).$
Since $G'_\bet(\atil)=0$ for $|\alp|<|\bet|$ it follows that
$G$ and $H$ are block triangular. Thus $h_{\bet\gam}=0$ for 
$|\gam|>|\bet|$ and we deduce that $b_{\bet\alp}=b_\bet(\abar^{-1})$
where $b_\bet:= \sum_{\gam:|\gam|\le|\bet|} h_{\bet\gam} O_\gam$
is a polynomial of degree $\le |\bet|$.

The top degree term on the left side of $(*)$ is a multiple
of $E_\alp$, and so by the definition of $G'_\alp$ we obtain that
$b_{\bet\alp}=0$ for $|\alp|\le|\bet|, \alp\ne\bet$. Thus
$b_\bet(\abar^{-1})=0$ for $|\alp|\le|\bet|, \alp\ne\bet$,
and since $\abar^{-1}=\abar(q^{-1},t^{-1})$, it follows that 
$b_\bet(x)$ is a multiple of $G_\bet(x;q^{-1},t^{-1}).$ In other
words, there are scalars $c_\bet$ in $\KK$ such that
$$ 
\qquad {G_\alp(ax)\over G_\alp(a\tau)}
= \sum_\bet c_\bet \qtinv
G'_\bet(x).$$

Comparing the top degree terms we get 
$c_\alp=a^{|\alp|}/G_\alp(a\tau)$ and the result follows. \qed

\Proof: (of \cite{binom2}) The proof proceeds similarly using
\cite{Oko2}. \qed

\beginsection Jack. More on the Jack limit

We now prove \cite{relate} and the symmetric versions of \cite{binom2} 
and \cite{relate}. Since the $(q,t)$-case will not
be considered in this section, we will often suppress $r$ 
to simplify the notation, e.g. we will write $G_\alp(x)$ 
for $G_\alp(x;r)$, $\btil$ for $\btil(r)$ etc.

We start with a simple, but crucial, lemma.

\Lemma dom. $1)$ $w_{-w_o\bet}=w_ow_\bet w_o$; $2)$ $-w_o\btil= 
\bbar+(n-1)r$.

\Proof: For $w$ in $S_n$, we have  $(w_oww_o)^{-1}(-w_o\bet)=
(-w_o)(w^{-1}\bet)$, which is dominant if and only if 
$w^{-1}\bet$ is dominant. Since conjugation by $w_o$
preserves length, part 1) follows. 

Now $\btil=\overline{-w_o\bet}=
-w_o\bet+w_{-w_o\bet}\rho= -w_o\bet+w_ow_\bet w_o\rho$, by part 1).
Also, since  $w_o\rho=-(n-1)r-\rho$, we get
$\btil=-w_o( \bet+(n-1)r+ w_\bet\rho)=-w_o( \bbar+(n-1)r)$.\qed

\Proof: (of \cite{relate}) For any polynomial $f$, $wf$ and 
$\sig(w)f$ have the same top terms. So, since $w_0^2=1$, 
the top term on the right of \cite{relate} is 
$(-1)^{|\alp|}w_0^2E_\alp(-x)=E_\alp(x)$.
It remains only to show that the right side of \cite{relate} 
belongs to the space $V$ consisting of polynomials which vanish 
at the points $x=\btil, |\bet|<|\alp|$. 

Putting $a=0$ and $v=\bet$ in 2) of \cite{recur2}, we 
deduce that $V$ is $\sig$-invariant and so it suffices
to prove that $f\equiv w_oG_\alp\left(-x-(n-1)r\right)\in V$. 
But, using \cite{dom} we get $$f(\btil)= w_oG_\alp
\left(-\btil-(n-1)r\right)=G_\alp\left(-w_o\btil-
(n-1)r\right)=G_\alp(\bbar),$$ which vanishes
for $|\bet|<|\alp|$ by the definition of $G_\alp$.  \qed

We now turn to the symmetric versions of \cite{binom2} 
and \cite{relate}.  As in \cite{OO} we define the ``symmetric 
$r$-binomial 
coefficients'' by 
$$ \lmr={R_\mu(\lbar(r);r)\over R_\mu(\mbar(r);r)}.$$

The main result of \cite{OO} is the generalized binomial formula
$$(**) \qquad {P_\lam(1+x;r)\over P_\lam(1;r)}
= \sum_{\mu\subseteq\lam} \lmr 
{P_\mu(x;r)\over P_\mu(1;r)}.$$

For the inhomogeneous analogue of this result, we define

\Definition: $R'_\lam(x;r)$ is the 
unique symmetric polynomial in $\QQ(r)[x]$ such that
\item{1)} $R'_\lam(x;r)$ and $R_\lam(x;r)$ have the same top degree terms,
\item{2)} $R'_\lam(x;r)$ vanishes at $x=\mtil(r)\equiv \overline{-w_o\mu}(r)$ 
for all $\mu$ with $|\mu|<|\lam|$.

Then we have 

\Theorem relate2. $\displaystyle R'_\lam(x;r)=
(-1)^{|\lam|}R_\lam(-x-(n-1)r;r).$

\Proof: The two sides have the same top degree terms,
and it suffices to prove that the right side vanishes 
for $x=\mtil$ if $|\mu|<|\lam|$.
By symmetry, we may consider instead $x=w_o\mtil$.
Substituting this and using \cite{dom}, the
right side becomes $(-1)^{|\lam|}R_\lam(\mbar;r)$,
which vanishes by definition of $R_\lam$. \qed

\Theorem binom3. $\displaystyle {R_\lam(a+x;r)\over 
R_\lam(a+\rho;r)} = \sum_{\mu\subseteq\lam} \lmr 
{R'_\mu(x;r)\over R_\mu(a+\rho;r)}.$

We shall deduce \cite{binom3} from \cite{binom2} by symmetrization.
Write $\cS$ for the operator
${1\over n!}\sum_{w\in S_n}\sig(w)$ acting on $\QQ(r)[x]$.

\Lemma sym. $\cS$ maps polynomials to symmetric polynomials. 

\Proof: For all $i$, we have $\sig_i\cS=\sum_{w\in S_n}\sig(s_iw)=\cS$.
So if $f$ is a polynomial in the image of $\cS$, then $(1-\sig_i)f=0$.
Rewriting this we get $(1- {r\over x_i-x_{i+1}})(1-s_i)f=0$. Hence
$(1-s_i)f=0$ for all $i$, which implies that $f$ is symmetric. \qed

\Lemma symm. Let $\alp$ be any composition with $\alp^+=\lam$,
then 
$$1)\; \; 
{\cS G_\alp(a+x)\over G_\alp(a+\rho)}= 
{R_\lam(a+x)\over R_\lam(a+\rho)}; 
\quad  2)\; \; {\cS G'_\alp(x)\over G_\alp(a+\rho)}= 
{R'_\lam(x)\over R_\lam(a+\rho)}.$$ 

\Proof: If $|\bet|\le |\alp|$ and $\bet^+\ne \alp^+$, then \cite{discr2} 
implies that, for all $w$ in $S_n$, the polynomial $\sig(w)G_\alp(x)$ 
vanishes at $x=\bbar$. This means that $f=\cS G_\alp(a+x)$ vanishes at 
$\mbar-a$ for all partitions $\mu$ satisfying $|\mu|\le|\lam|, 
\mu\ne\lam$. Since $f$ is symmetric and of the right degree, we 
conclude that $f$ is a multiple of $R_\lam(a+x)$. To determine the 
multiple we merely evaluate both sides of 1) at $x=\rho$, and use the fact 
that $\sig(w)G_\alp(a+\rho)=G_\alp(a+\rho)$ which follows from
2) of \cite{discr2}. This proves 1).

For 2), the same argument proves that $\cS G'_\alp(x)$ vanishes at 
$\mtil$ for $|\mu|<|\lam|$. To finish the proof, it suffices 
to prove that the {\it top\/} terms of the two sides are equal.
But these are also the top terms of 1) and hence are 
equal. \qed

\Proof: (of \cite{binom3}) Fix $\alp$ with $\alp^+=\lam$ and
apply $\cS$ to both sides of \cite{binom2}. By  \cite{symm} we get 
$$\displaystyle {R_\lam(a+x;r)\over R_\lam(a+\rho;r)}
= \sum_{\mu\subseteq\lam} k_\mu 
{R'_\mu(x;r)\over R_\mu(a+\rho;r)}\; \|with| 
\; k_\mu=\sum_{\bet^+=\mu} \abr\in \QQ(r).$$

To conclude we need to establish that $k_\mu=\lmr$, but this 
follows by putting $x=ax$ in the above, letting $a\rightarrow\infty$,
and using $(**)$. \qed

\Corollary rel. For each $\alp$ satisfying $\alp^+=\lam$, 
we have $\displaystyle \sum_{\bet^+=\mu} \abr=\lmr. $ \qed

\beginsection References. References

\baselineskip12pt
\parskip2.5pt plus 1pt
\hyphenation{Hei-del-berg}
\def\L|Abk:#1|Sig:#2|Au:#3|Tit:#4|Zs:#5|Bd:#6|S:#7|J:#8||{%
\edef\TEST{[#2]}
\expandafter\ifx\csname#1\endcsname\TEST\relax\else
\immediate\write16{#1 had already been defined!}\fi
\expandwrite\AUX{\neverexpand\ref{#1}{\TEST}}
\HI{[#2]}
\ifx-#3\relax\else{#3}: \fi
\ifx-#4\relax\else{#4}{\sfcode`.=3000.} \fi
\ifx-#5\relax\else{\it #5\/} \fi
\ifx-#6\relax\else{\bf #6} \fi
\ifx-#8\relax\else({#8})\fi
\ifx-#7\relax\else, {#7}\fi\Par}

\def\B|Abk:#1|Sig:#2|Au:#3|Tit:#4|Reihe:#5|Verlag:#6|Ort:#7|J:#8||{%
\edef\TEST{[#2]}
\expandafter\ifx\csname#1\endcsname\TEST\relax\else
\immediate\write16{#1 hat sich geaendert!}\fi
\expandwrite\AUX{\neverexpand\ref{#1}{\TEST}}
\HI{[#2]}
\ifx-#3\relax\else{#3}: \fi
\ifx-#4\relax\else{#4}{\sfcode`.=3000.} \fi
\ifx-#5\relax\else{(#5)} \fi
\ifx-#7\relax\else{#7:} \fi
\ifx-#6\relax\else{#6}\fi
\ifx-#8\relax\else{ #8}\fi\Par}

\def\Pr|Abk:#1|Sig:#2|Au:#3|Artikel:#4|Titel:#5|Hgr:#6|Reihe:{%
\edef\TEST{[#2]}
\expandafter\ifx\csname#1\endcsname\TEST\relax\else
\immediate\write16{#1 hat sich geaendert!}\fi
\expandwrite\AUX{\neverexpand\ref{#1}{\TEST}}
\HI{[#2]}
\ifx-#3\relax\else{#3}: \fi
\ifx-#4\relax\else{#4}{\sfcode`.=3000.} \fi
\ifx-#5\relax\else{In: \it #5}. \fi
\ifx-#6\relax\else{(#6)} \fi\PrII}
\def\PrII#1|Bd:#2|Verlag:#3|Ort:#4|S:#5|J:#6||{%
\ifx-#1\relax\else{#1} \fi
\ifx-#2\relax\else{\bf #2}, \fi
\ifx-#4\relax\else{#4:} \fi
\ifx-#3\relax\else{#3} \fi
\ifx-#6\relax\else{#6}\fi
\ifx-#5\relax\else{, #5}\fi\Par}
\setHI{[ABC]\ }
\sfcode`.=1000

\L|Abk:Bi|Sig:B|Au:Bingham, C.|Tit:An identity involving
partitional generalized binomial coefficients|Zs:J.~Multiv.
Anal.|Bd:4|S:210--223|J:1974||

\L|Abk:C|Sig:C|Au:Cherednik, I.|Tit:Nonsymmetric Macdonald
polynomials|Zs:IMRN|Bd:10|S:483--515|J:1995||

\L|Abk:K|Sig:K|Au:Knop, F.|Tit:Symmetric and nonsymmetric 
quantum Capelli polynomials|Zs:-|Bd:-|S:-|J:preprint||

\L|Abk:KS|Sig:KS|Au:Knop, F.; Sahi, S.|Tit:Difference
equations and symmetric polynomials defined by their 
zeros|Zs:IMRN|Bd:10|S:473--486|J:1996||

%\L|Abk:KS|Sig:KS|Au:Knop, F.; Sahi, S.|Tit:A recursion 
%and a combinatorial formula for Jack polynomials|Zs:Invent. 
%Math.|Bd:-|S:-|J:to appear||

\L|Abk:Lc|Sig:Lc|Au:Lasocux, A.|Tit:Classes de Chern d'un
produit tensoriel|Zs:C.~R.~Acad.~Sci.|Bd:286A|S:385--387|J:1978||

\L|Abk:Ls|Sig:Ls|Au:Lassalle, M.|Tit:Une formule de 
bin\^ome g\'en\'eralis\'ee pour les polyn\^omes de 
Jack|Zs:C.~R.~Acad.~Sci. Ser.1|Bd:310|S:253--256|J:1990||

\B|Abk:M1|Sig:M1|Au:Macdonald, I.|Tit:Symmetric Functions 
and Hall Polynomials|Reihe:2nd ed.|Verlag:Clarendon 
Press|Ort:Oxford|J:1995||

\L|Abk:M2|Sig:M2|Au:Macdonald, I.|Tit:Affine Hecke algebras and 
orthogonal polynomials|Zs:S\'em. Bourbaki|Bd:797|S:1--18|J:1995||

\L|Abk:O|Sig:O|Au:Okounkov, A.|Tit:Binomial formula for 
Macdonald polynomials|Zs:-|Bd:-|S:-|J:preprint||

\L|Abk:OO|Sig:OO|Au:Okounkov, A., Olshanski, G.|Tit:Shifted 
Jack polynomials, binomial formula, and 
applications|Zs:-|Bd:-|S:-|J:preprint||

\L|Abk:Op|Sig:Op|Au:Opdam, E.|Tit:Harmonic analysis 
for certain representations of graded Hecke 
algebras|Zs:Acta Math.|Bd:175|S:75--121|J:1995||

\Pr|Abk:S1|Sig:S1|Au:Sahi, S.|Artikel:The spectrum of certain 
invariant differential operators associated to a Hermitian 
symmetric space|Titel:Lie theory and Geometry|Hgr:-|Reihe:Progr. 
Math|Bd:123|Verlag:Birkhauser|Ort:Boston|S:569--576|J:1994||

\L|Abk:S2|Sig:S2|Au:Sahi, S.|Tit:Interpolation, 
integrality, and a generalization of Macdonald's 
polynomials|Zs:IMRN|Bd:10|S:457--471|J:1996||

\L|Abk:S3|Sig:S3|Au:Sahi, S.|Tit:A new scalar product for 
nonsymmetric Jack polynomials|Zs:IMRN|Bd:20|S:997--1004|J:1996||

\bye